\documentclass[conference]{IEEEtran}
\IEEEoverridecommandlockouts
\usepackage{cite}
\usepackage{amsmath,amssymb,amsfonts,mathtools}
\usepackage{algorithm}
\usepackage[noend]{algpseudocode}
\usepackage{graphicx}
\usepackage{textcomp}
\usepackage{xcolor}
\def\BibTeX{{\rm B\kern-.05em{\sc i\kern-.025em b}\kern-.08em
    T\kern-.1667em\lower.7ex\hbox{E}\kern-.125emX}}

\DeclareMathOperator*{\argmax}{arg\,max}
\DeclareMathOperator*{\argmin}{arg\,min}
\DeclarePairedDelimiter\abs{\lvert}{\rvert}%

\begin{document}

\title{Monetizing Customer Load Data for an Energy Retailer: A Cooperative Game Approach \\
\thanks{This work was accepted and presented at the 14th IEEE PowerTech 2021.}
\thanks{The research leading to this work is being carried out as a part of the Smart4RES project (European Union’s Horizon 2020, No. 864337). The
sole responsibility of this publication lies with the author. The European
Union is not responsible for any use that may be made of the information
contained therein.}
}

\author{\IEEEauthorblockN{Liyang Han\IEEEauthorrefmark{1}, Jalal Kazempour\IEEEauthorrefmark{1}, Pierre Pinson\IEEEauthorrefmark{2}}
\IEEEauthorblockA{\textit{\IEEEauthorrefmark{1}Department of Electrical Engineering, \IEEEauthorrefmark{2}Department of Technology, Management and Economics} \\
\textit{Technical University of Denmark}\\
Kongens Lyngby, Denmark \\
liyha@elektro.dtu.dk, seykaz@elektro.dtu.dk, ppin@dtu.dk}
}

\maketitle

\begin{abstract}

When energy customers schedule loads ahead of time, this information, if acquired by their energy retailer, can improve the retailer's load forecasts. Better forecasts lead to wholesale purchase decisions that are likely to result in lower energy imbalance costs, and thus higher profits for the retailer. Therefore, this paper monetizes the value of the customer schedulable load data by quantifying the retailer's profit gain from adjusting the wholesale purchase based on such data. Using a cooperative game theoretic approach, the retailer translates their increased profit in expectation into the value of cooperation, and redistributes a portion of it among the customers as monetary incentives for them to continue providing their load data. Through case studies, this paper demonstrates the significance of the additional profit for the retailer from using the proposed framework, and evaluates the long-term monetary benefits to the customers based on different payoff allocation methods.

\end{abstract}

\begin{IEEEkeywords}
data monetization, energy market, cooperative game theory, newsvendor model, probabilistic forecasting
\end{IEEEkeywords}

\section{Introduction}

The increased penetration of renewable energy generation and the sustained deployment of distributed energy resources (electric vehicles, rooftop PV, batteries, etc.) are inducing profound changes to the operation and economics of electric energy systems. Keeping pace with these changes, digitization is unlocking access to a wealth of data for improved observability, forecasting, data-driven decision making and control. In order to minimize the disruption to the operation of energy systems, adequate governance and business models are required to ensure that all relevant agents are on board with the green transition \cite{LeRay2020}.

Focusing on the market, energy retailers make a profit by purchasing energy in the wholesale market and reselling it to their customers at a higher price. They are subject to imbalance costs, a form of financial penalty for the mismatch between a retailer's advance purchase of wholesale energy and their customers' real-time demand \cite{EuropeanCommission2017}. It is therefore in the retailer's best interest to improve their forecast of the demand in order to determine the optimal purchase amount of wholesale energy that maximizes their expected profit \cite{Carrion2009}. Meanwhile, the newly gained autonomy from distributed generation \cite{Scheer2012} and controllability from smart appliances and energy storage systems are fueling the energy customers' interest in becoming active participants in the market operation \cite{Parag2016}. Recent literature on this topic heavily focuses on the creation of distributed energy markets, which generally allow customers to share or trade energy among themselves \cite{Sousa2019}. Recognizing that the way customers control and schedule their energy resources holds important information that can help energy retailers to improve their demand forecast, we identify customer load data as another commodity that creates value, and thus has the potential to be monetized in the energy market. 

A common method of collecting energy customer load data is through the use of smart meters. Privacy is a major concern despite research effort in developing privacy-preserving techniques to process customer data \cite{Farokhi2020}. The lack of fairness is another important issue considering the value of data is completely absorbed by the data collector, who is the retailer in an energy market context. To address these issues, we adopt the concept of a \emph{data market}, in which customers are financially compensated for disclosing their data \cite{Bergemann2019}. Specifically, we consider the retailer a \emph{data buyer} and the customers \emph{data sellers}, and we focus on the monetization of customer load data and the allocation of payoffs. In a real-world application, customers can choose to reject the financial offer when it is lower than their privacy value, but the valuations of privacy, and of privacy leakage as a result of data correlation among data sellers \cite{Acemoglu2019}, are beyond the scope of this paper. 

To determine the data buyer's financial offers in a general data market framework, an auction can be implemented to achieve market equilibrium based on noncooperative game theory \cite{Bimpikis2019}. Alternatively, a cooperative game can be adopted to determine the allocation of payoffs based on the contribution of each seller's data to the joint welfare \cite{Agarwal2019}. Using this cooperative game mechanism,  a data platform is constructed in \cite{Goncalves2020} for wind generators to trade power measurement data with each other in order to improve their local forecasts of wind power outputs. Nevertheless, the total payment as a portion of the data buyer's profit has to be predetermined by the market operator in favor of the buyer to ensure no regret, undermining fairness, which is a desirable outcome of a cooperative game.

In order to incentivize energy customers to offer their load data to help the energy retailer optimize their wholesale purchase, hence increasing their profit, we propose:
\begin{itemize}
\item A probabilistic customer load demand forecast based on the customers' schedulable load data;
\item An adapted \emph{newsvendor model} \cite{Edgeworth1888} to represent an energy retailer, who uses their probabilistic forecast of the total customer demand to determine the amount of wholesale purchase that maximizes their profit in expectation;
\item A data market built upon a cooperative game framework to monetize customer load data. Different from the framework proposed in \cite{Agarwal2019}, which predetermines the total data seller profit to ensure no regret for the data buyer, we include the data buyer, the retailer, in the game in order to achieve full cooperation.
\end{itemize}

The structure of the paper is as follows. Section \ref{sec:ener_newsv} sets up a cost-based newsvendor model for the energy retailer and establishes the link between customer schedulable load data and the retailer's probabilistic forecast of the total demand. Then, Section \ref{sec:coop} introduces a cooperative game approach to monetize customer schedulable load data and details the payoff allocations. Section \ref{sec:case_studies} presents several case studies to demonstrate the effectiveness of the proposed framework in increasing the retailer's profit while incentivizing the customers to provide their schedulable load data to the retailer. Finally, Section \ref{sec:conc} concludes the paper.

\section{Energy Retailer Newsvendor Model} \label{sec:ener_newsv}

We consider a monopolistic energy retailer model, in which the energy retailer makes advance purchases of energy in the wholesale  market and then resells the energy to the retail customers in real-time at a higher price. It is shown in this section that this energy retailer is analogous to a cost-based newsvendor, and by choosing the wholesale purchase based on the retailer's probabilistic forecast of the customer loads, the retailer is able to maximize their profit in expectation.

\subsection{Energy Retailer Model Setup}

Similar to the retail business models of commercial goods, the profit of an energy retailer comes from the mark-up of the energy price (i.e., retail price $>$ wholesale price). However, one main difference between an energy retailer and a retailer of other commercial goods is the real-time operational constraint of an energy network, where supply has to be always equal to demand. Therefore, even though the energy retailer is free to choose the amount of wholesale energy to purchase, any surplus or shortage of wholesale purchase compared to the real-time demand has to be absorbed or replenished by the grid's regulation services. The costs of these services are then passed on to the retailer in the form of \emph{imbalance costs}.

Even though the energy network operation is continuous in time, the energy market is cleared in discrete time units (e.g., 30-minute, 15-minute, etc.). For each time unit, we assume the retailer makes an advance purchase of $q$ amount of energy (kWh) in the wholesale market based on the forecast of the customers' total demand. The energy wholesale market usually consists of multiple stages (e.g., day-ahead, intraday, hour-ahead, etc.). The retailer's bids in these stages all depend on the customer load forecasts, so we combine them into one stage and use the average price $r^\mathrm{w}$ as the wholesale price for simplicity. The retailer then charges the customers the retail price $r^\mathrm{r}$ ($r^\mathrm{r} > r^\mathrm{w}$) for their real-time demand $D$. Based on the mismatch between $q$ and $D$, an imbalance cost is imposed on the retailer. The retailer's profit $H$ for a single time unit is
\begin{align}  \label{eq:enerv_obj}
H(q, D) = \ & r^\mathrm{r} D - r^\mathrm{w} q  - b^- \max(D-q, 0)   \nonumber \\
 & + b^+ \max(q-D, 0),
\end{align}
where $b^-$ and $b^+$ are the imbalance prices for a shortage and a surplus, respectively, of the purchased wholesale energy compared to the real-time demand. Imbalance prices can be either positive or negative, and we define their signs to be consistent with the European Commission (see Table \ref{tb:imbalance_price}) \cite{EuropeanCommission2017}.

\begin{table}[t]
\begin{center}
\caption{Payment for Balancing Energy}  \label{tb:imbalance_price}
\begin{tabular}{c|c|c}
\hline
& $>0$ &$<0$ \\\hline
$b^+$ & Retailer receives payment & Retailer makes payment \\\hline
$b^-$ & Retailer makes payment &  Retailer receives payment \\\hline
\end{tabular}
\end{center}
\end{table}

The retail demand $D$ is a random variable following a probability density function (PDF) $f_D$, and a corresponding cumulative distribution function (CDF) $F_D(d) = P(D \leq d)$. Given the retailer model is probabilistic, the objective of the retailer becomes to maximize their profit in expectation, i.e., $\max_q \mathbb{E}[H(q, D)]$, where $\mathbb{E}$ is the expectation operator.

We notice that the revenue term $r^\mathrm{r} D$ in \eqref{eq:enerv_obj} is not affected by the wholesale purchase $q$, so we can focus on only the retailer's costs, i.e.,
\begin{equation}
C(q, D) = r^\mathrm{w} q + b^- \max(D-q, 0) - b^+ \max(q-D, 0),
\end{equation}
and reconstruct the retailer's profit maximization problem as an equivalent cost minimization problem, i.e.,
\begin{align} \label{eq:enerv_cost}
\min_q \mathbb{E}[C(q, D)] =  \min_q \{ r^\mathrm{w} q & + b^- \mathbb{E}[\max(D-q, 0)] \nonumber \\
& - b^+ \mathbb{E}[\max(q-D, 0)]\}.
\end{align}

The format of this energy retailer model is the same as the cost-based newsvendor model discussed in \cite{Birge2011}, following which we can obtain the solution for \eqref{eq:enerv_cost} as
\begin{equation} \label{eq:enerv_opt_quan}
q^\ast = \argmin_q \mathbb{E}[C(q, D)] = F_D^{-1} \left( \gamma \right),
\end{equation}
where $\gamma = \frac{b^- - r^\mathrm{w}}{b^- - b^+}$ is the retailer's \emph{cost ratio} based on Littlewood's rule \cite{Littlewood2005}.

Imbalance prices are unknown a priori as they are usually determined based on the real-time imbalance of the energy network. Forecasting of the imbalance prices is needed in practice to calculate the cost ratio, but it is not the focus of this paper. To simulate the proposed framework, we adopt a dual pricing scheme ($b^- \neq b^+$) \cite{Pinson2007} and assume known imbalance prices that satisfy $b^- \geq r^\mathrm{w}$ and $b^+ \leq r^\mathrm{w}$. The rationale behind these inequalities is that the production of balancing generators is more costly than the successful bids in the wholesale market, and that those with surplus energy would be encouraged to trade out their positive imbalances in advance \cite{Morales2014}.

\subsection{Probabilistic Customer Loads}

Based on the energy retailer newsvendor model, the optimal wholesale purchase $q^\ast$ is dependent on $F_D(d)$, hence the probability distribution of the total customer demand $D$. 

Assume the energy retailer has $M$ customers, the set of which is denoted as $\mathcal{M}$, and each customer $i$ $(i \in \mathcal{M})$ has two types of loads: schedulable loads $L^\mathrm{s}_i$ and unschedulable loads $L^\mathrm{u}_i$. Schedulable loads can be determined at the customer's discretion, unbeknownst to the retailer, before the gate closure of the wholesale market, whereas the unschedulable loads are unknown to both the customers and the retailer. If the retailer's probabilistic forecasts of the schedulable and unschedulable loads are described by PDFs $f_{L^\mathrm{s}_i}$ and $f_{L^\mathrm{u}_i}$ respectively, the PDF and the CDF of the total demand D are
\begin{align}
f_D =& f_{\sum_{i \in \mathcal{M}} (L^\mathrm{s}_i+L^\mathrm{u}_i)}, \\
F_D =& F_{\sum_{i \in \mathcal{M}} (L^\mathrm{s}_i+L^\mathrm{u}_i)}.
\end{align}

In a special case where every customer's schedulable and unschedulable loads follow the Gaussian distribution as $L^\mathrm{s}_i \sim \mathcal{N}(\mu^\mathrm{s}_i, {\sigma^\mathrm{s}_i}^2)$ and $L^\mathrm{u}_{i} \sim \mathcal{N}(\mu^\mathrm{u}_i, {\sigma^\mathrm{u}_i}^2)$, the retailer's forecast of the total customer load follows
\begin{equation} \label{eq:data_emptyset}
D_\emptyset \sim  \  \mathcal{N} \left( \sum_{i \in \mathcal{M}} (\mu^\mathrm{s}_i + \mu^\mathrm{u}_i), \sum_{i \in \mathcal{M}}({\sigma^\mathrm{s}_i}^2 + {\sigma^\mathrm{u}_i}^2) \right),
\end{equation}
where the subscript $\emptyset$ indicates that none of the customers have disclosed their schedulable load data to the retailer. Therefore, \eqref{eq:enerv_opt_quan} can be solved analytically as
\begin{equation}  \label{eq:opt_proc_emptyset}
q^\ast_\emptyset = \Phi^{-1} (\gamma) \sqrt{\sum_{i \in \mathcal{M}} ({\sigma^\mathrm{s}_i}^2 + {\sigma^\mathrm{u}_i}^2)}  +  \sum_{i \in \mathcal{M}} (\mu^\mathrm{s}_i + \mu^\mathrm{u}_i),
\end{equation}
where $\Phi^{-1}$ is the inverse of the standard Gaussian CDF. 

\subsection{Sharing of Schedulable Load Data} \label{subs:data_sharing}

Since the retailer's profit is influenced by the accuracy of their forecast of the total customer demand, there is an incentive for the retailer to obtain the data on how the customers schedule their loads. Assuming a subset of customers $\mathcal{T}$ $(\mathcal{T} \subseteq \mathcal{M})$ are willing to disclose to the retailer the data on their schedulable loads $l^\mathrm{s}_i, \forall i \in \mathcal{T}$, their schedulable loads become deterministic in the retailer's forecast. We define $\mu_\mathcal{T}$ and $\sigma_\mathcal{T}^2 \ (\sigma_\mathcal{T} \geq 0)$ as the mean and the variance of the retailer's forecast of the total demand given $\mathcal{T}$'s schedulable load data, so we have
\begin{align}
\mu_\mathcal{T} = & \sum_{i \in \mathcal{T}} l^\mathrm{s}_i + \sum_{i \in (\mathcal{M} \setminus \mathcal{T})} \mu^\mathrm{s}_i + \sum_{i \in \mathcal{M}} \mu^\mathrm{u}_i, \\
\sigma_\mathcal{T}^2 = & \sum_{i \in (\mathcal{M} \setminus \mathcal{T})} {\sigma^\mathrm{s}_i}^2 + \sum_{i \in \mathcal{M}} {\sigma^\mathrm{u}_i}^2.
\end{align}
To be consistent with \eqref{eq:data_emptyset} and \eqref{eq:opt_proc_emptyset}, we denote $D_\mathcal{T}$ as the retailer's updated forecast given the schedulable load data from $\mathcal{T}$, and $q^\ast_\mathcal{T}$ as the corresponding optimal wholesale purchase:
\begin{equation}
D_\mathcal{T} \sim    \mathcal{N} \left(\mu_\mathcal{T}, \sigma_\mathcal{T}^2\right),
\end{equation}
\begin{equation}
q^\ast_\mathcal{T} = \Phi^{-1} (\gamma) \sigma_\mathcal{T} + \mu_\mathcal{T}.
\label{eq:opt_bid_cl}
\end{equation}

As an example, we randomly select 8 customer load profiles from the domestic load data obtained in the UK Customer-Led Network Revolution trials\footnote{http://www.networkrevolution.co.uk/project-data-download/?dl=TC5.zip}, so $\mathcal{M} = \{1,...,8\}$. For each timestep, we implement the following:
\begin{enumerate}
\item The measured load $\mu_i$ of each customer $(i \in \mathcal{M})$ is used as the mean in the retailer's forecast if no schedulable load data are disclosed to the retailer: $\mu_\emptyset = \sum_{i \in \mathcal{M}} \mu_i$. 
\item $\mu_i$ is separated into a schedulable part and an unschedulable part: $\mu_i = \mu^\mathrm{s}_i + \mu^\mathrm{u}_i$, where the schedulable load is determined as a fraction of $\mu_i$, i.e., $\mu^\mathrm{s}_i=\alpha_{si} \mu_i$. 
\item The standard deviations of the schedulable and unschedulable loads are determined as fractions of the respective loads: ${\sigma^\mathrm{s}_i} = \beta_{si} \mu^\mathrm{s}_i, \ \sigma^\mathrm{u}_{i} = \beta_{ui} \mu^\mathrm{u}_i$. 
\item Without any schedulable load data, the retailer's forecast of the total load is $D_ \emptyset \sim \mathcal{N} \big(\mu_\emptyset, \sum_{i \in \mathcal{M}}({\sigma^\mathrm{s}_i}^2 + {\sigma^\mathrm{u}_i}^2) \big)$.
\item Customers determine their schedulable loads $l^\mathrm{s}_i$, which are realizations of $L^\mathrm{s}_i \sim \mathcal{N}(\mu^\mathrm{s}_i, {\sigma^\mathrm{s}_i}^2), \forall i \in \mathcal{M}$. Given all customers' schedulable load data, the mean of the retailer's forecast becomes $\mu_\mathcal{M} = \sum_{i \in \mathcal{M}}(l^\mathrm{s}_i + \mu^\mathrm{u}_i)$.
\item With all the customers' schedulable load data, the retailer's forecast of the total load is $D_ \mathcal{M} \sim \mathcal{N}(\mu_\mathcal{M}, \sum_{i \in \mathcal{M}} {\sigma^\mathrm{u}_i}^2)$.
\end{enumerate}

\begin{figure} 
\begin{center}
\includegraphics[width=8.4cm]{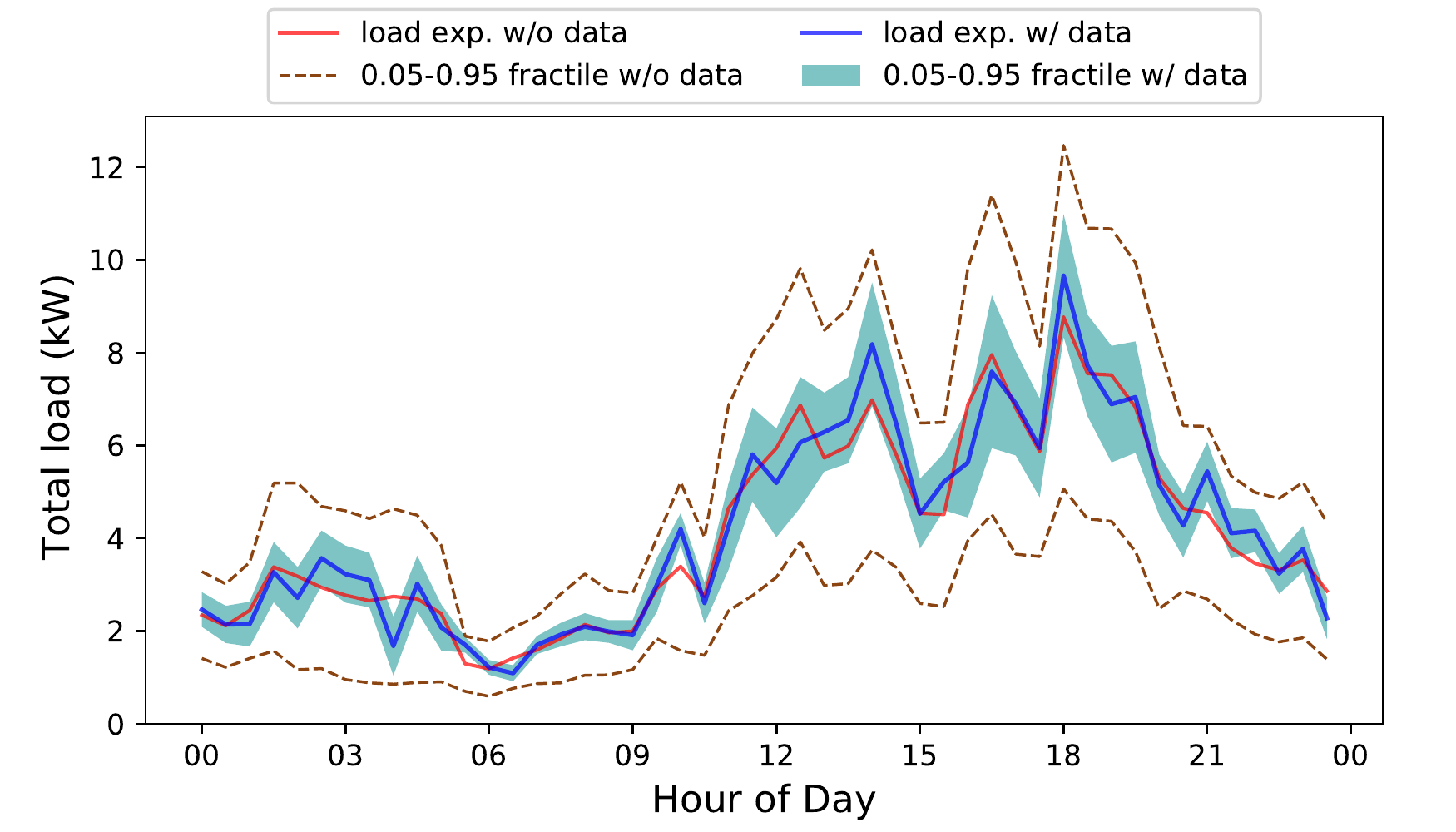}    
\caption{Retailer's forecasts of a 24-hour combined load profile of 8 customers: with and without customers' shared data on their schedulable loads.} 
\label{fig:l_prof}
\end{center}
\end{figure}

If we assign $\beta_{si} = 1.0,  \beta_{ui} = 0.5$, and randomly select values between 10-90\% for $\alpha_{si}, \forall i \in \mathcal{M}$ for each timestep, we obtain the retailer's forecast of the total load with and without data, as shown in Fig. \ref{fig:l_prof} for a January day. The mean profiles of the two forecasts are slightly different because the retailer with data no longer needs to forecast the schedulable loads. Additionally, the same confidence interval covers a much smaller range of values for the retailer with schedulable load data. This difference in forecasting leads to different wholesale purchase decisions hence different profit expectations, which are further investigated in Section \ref{sec:coop}.

\section{Proposed Cooperative Game Formulation} \label{sec:coop}

Schedulable load data improve the retailer's forecast, which leads to purchase decisions that result in a higher expected profit. However, customers are unlikely to offer their data for free. In order to incentivize the customers to provide their schedulable load data, the retailer can offer a share of their additional profit in return. For fairness, the financial rewards should reflect the added value of the data provided by each customer. To this end, we adopt cooperative game theory to evaluate the contributions of each customer's schedulable load data, and to allocate the profit based on these contributions.

A standard cooperative game consists of two main steps: 1) formation of coalitions and calculation of coalition values; 2) allocation of payoffs based on coalition values. These two steps are explained in detail in the following.

\subsection{Coalition Definition and Coalition Value Function}

Keeping the notations from Section \ref{sec:ener_newsv}, we define the set of all $M$ customers as $\mathcal{M} = \{1, ..., M\}$, where each customer is considered a player in the game. Including the retailer as `Player $0$', we define the \emph{grand coalition} of this cooperative game as $\mathcal{M}^0 = \{0, 1,...,M\}$. A coalition $\mathcal{T}$ is defined as any subset of the grand coalition: $\mathcal{T} \subseteq \mathcal{M}^0$. 

The value of cooperation in this proposed game is achieved through the customers sharing data with the retailer, and the retailer making more informed decisions on their wholesale purchases. Therefore, the participation of the retailer is essential for the value creation. In other words, coalitions without the retailer have zero values. On the other hand, if the retailer is part of a coalition $\mathcal{T}$ ($0 \in \mathcal{T}$), based on \eqref{eq:opt_bid_cl}, the optimal wholesale purchase is updated from $q^\ast_\emptyset$ to $q^\ast_{\mathcal{T} \setminus \{0\}}$. The profit gain for having customers in $\mathcal{T} \setminus \{0\}$ sharing their schedulable data with the retailer is defined as the value of coalition $\mathcal{T}$. These two scenarios are summarized as 
\begin{equation}
v(\mathcal{T}) = 
	\begin{cases} 
	0 ,&  \text{if} \ 0 \notin \mathcal{T}; \\
	\mathbb{E}[H(q^\ast_{\mathcal{T} \setminus \{0\}}, D_{\mathcal{M}})] - \mathbb{E}[H(q^\ast_\emptyset, D_{\mathcal{M}})]  ,&  \text{if} \ 0 \in \mathcal{T}. \\
	\end{cases}
\nonumber
\end{equation}
Notice that the total customer demand is realized based on the the scenario where all schedulable load data are known, i.e., 
\begin{equation} \label{eq:D_dist}
D_ \mathcal{M} \sim \mathcal{N} \left( \mu_\mathcal{M}, \sum_{i \in \mathcal{M}} {\sigma^\mathrm{u}_i}^2 \right).
\end{equation}

\begin{figure} 
\begin{center}
\includegraphics[width=8.4cm]{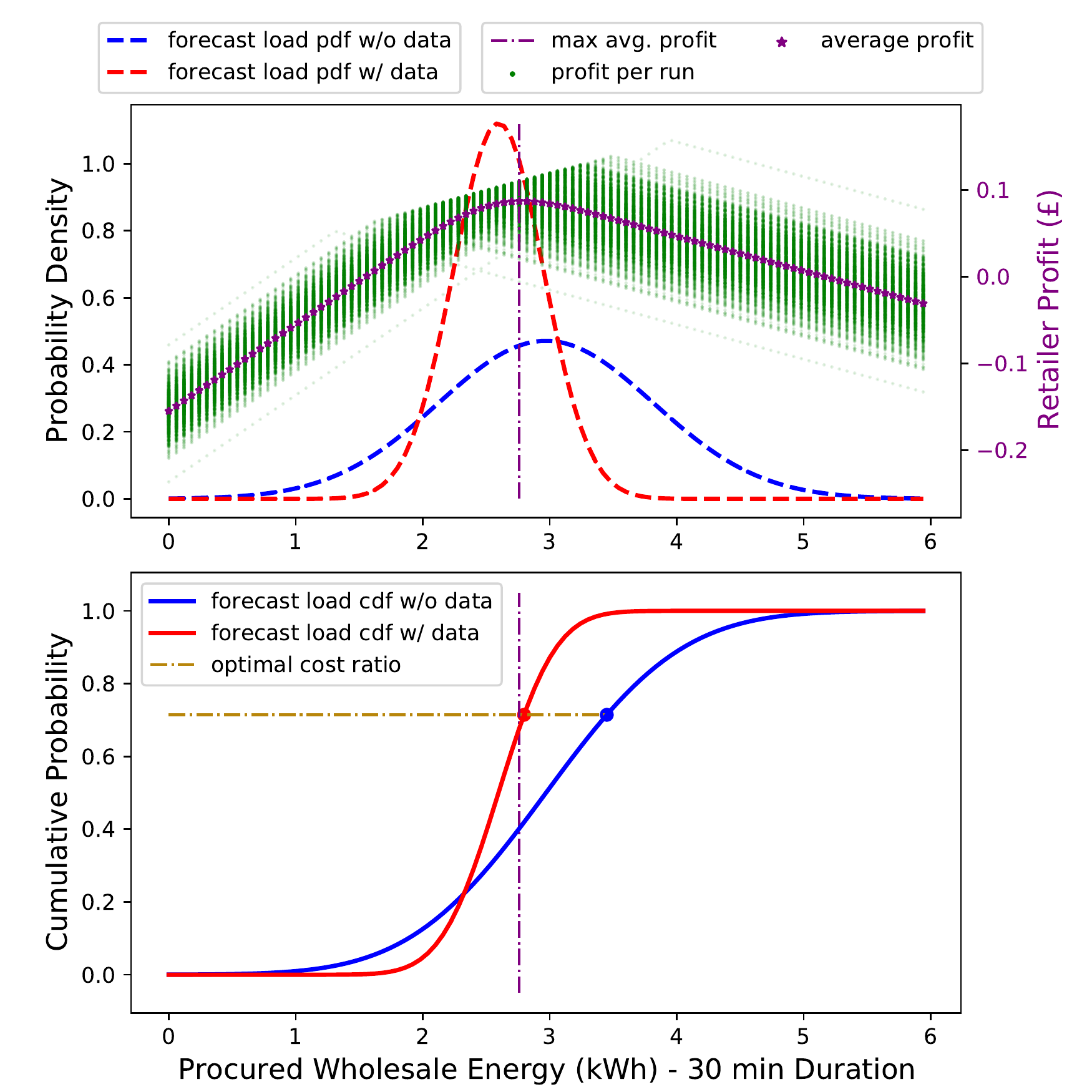}    
\caption{Single timestep snapshot of the retailer's load forecasts (with and without schedulable load data) and potential profits based on the procured wholesale energy.} 
\label{fig:pdf_cdf_prof}
\end{center}
\end{figure}

In this paper, we fix $r^\mathrm{r} = \pounds0.1$/kWh, $r^\mathrm{w} = \pounds0.06$/kWh, $b^- = \pounds0.16$/kWh, $b^+ = \pounds0.03$/kWh, and to calculate the expectation of profit $\mathbb{E}[H(q, D_ \mathcal{M})]$ given $q$, we take a Monte Carlo sampling approach with $D_ \mathcal{M}$ as the random variable. As shown in the top plot in Fig. \ref{fig:pdf_cdf_prof}, we show $[0, 2 \mu_\emptyset]$ as a reasonable range of purchased wholesale energy quantity, which is then evenly divided into 100 increments, each being a potential purchased wholesale quantity $q_j \ (j=1,...,100)$. For each $q_j$, $H(q_j, D_ \mathcal{M})$ is computed by sampling $D_ \mathcal{M}$ 1000 times following \eqref{eq:D_dist} , plotted as green dots. An average is taken and marked with a purple star, which represents $\mathbb{E}[H(q_j, D_ \mathcal{M})]$. We find $j^\ast = \argmax_{j} \mathbb{E}[H(q_j, D_ \mathcal{M})]$, which is the index for the quantity that contributes to the highest expected profit. $q_{j^\ast}$ is marked with the vertical purple line. On the other hand, we apply \eqref{eq:enerv_opt_quan} to identify the optimal wholesale quantity without data $q^\ast_\emptyset$ and with data $q^\ast_\mathcal{M}$, which are marked by the intersection of the brown horizontal line and the CDFs of each case in the bottom plot of Fig. \ref{fig:pdf_cdf_prof}. It is clear that the $q^\ast_\mathcal{M}$ coincides with $q_{j^\ast}$, proving the effectiveness of the Monte Carlo simulation approach. From this figure, we can also see that the retailer can improve their profit by updating their wholesale purchase from $q^\ast_\emptyset$ to $q^\ast_\mathcal{M}$.

\subsection{Profit Allocation} \label{}

To incentivize the disclosure of the schedulable load data, it is reasonable for the retailer to give back a portion of their profit gain to the customers. We use vector $\mathbf{x}$ as the \emph{payoff allocation}, and its entry $x_i$ represents the profit allocated to each player $i$. Based on the coalition value function definition, the total profit gain is $v(\mathcal{M}^0)$, so we require
\begin{equation}
\sum_{i \in \mathcal{M}^0} x_i = v(\mathcal{M}^0).
\end{equation}
This is the \emph{efficiency} criterion that guarantees that the entirety of the profit gained through cooperation gets allocated to the participants. We investigate two popular payoff allocations: the Shapley value and the nucleolus, both of which include the efficiency criterion in their computation. 

For the proposed game, the Shapley value is given by
\begin{equation}
\phi_i = \sum_{\mathclap{{\mathcal{T} \subseteq \mathcal{M}^0, i \in \mathcal{T}}}} \frac{(\abs{\mathcal{T}} - 1)! (\abs{\mathcal{M}^0} - \abs{\mathcal{T}})!}{\abs{\mathcal{M}^0}!} \left[v(\mathcal{T}) - v(\mathcal{T} \setminus \{i\})\right],
\label{eq:shapley_calc}
\end{equation}
where $\abs{\mathcal{M}^0} = M+1$, for $\mathcal{M}^0 = \{0, 1,...,M\}$.

The nucleolus is solved through an iteration process that minimizes the lexicographical coalition excesses \cite{Schmeidler1969}, a measure of coalitions' dissatisfaction given a payoff allocation. 

\begin{figure} 
\begin{center}
\includegraphics[width=8.4cm]{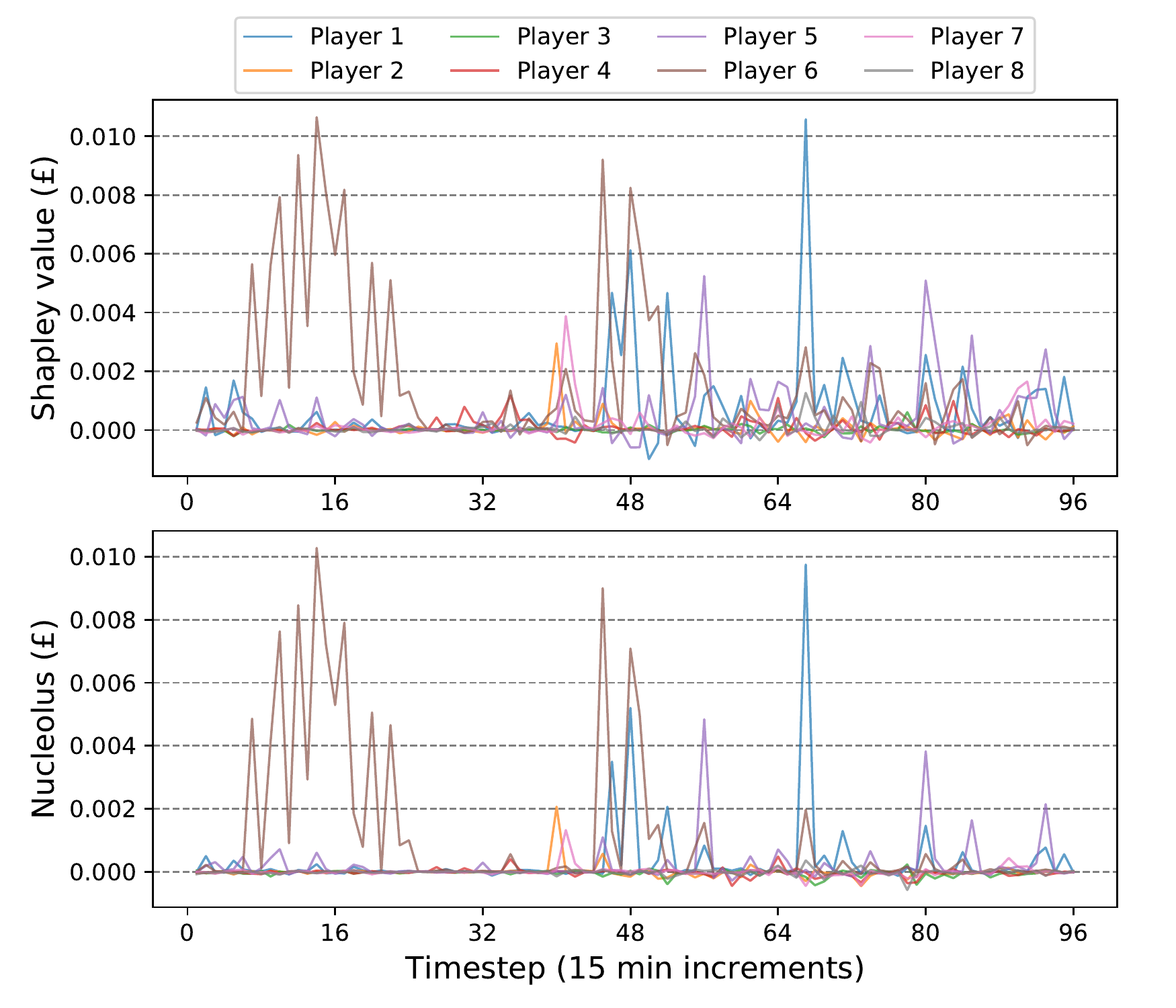}    
\caption{Half-hourly Shapley value and nucleolus allocated to all customers in a data sharing cooperative game over a 24-hr time span.} 
\label{fig:sh_t}
\end{center}
\end{figure}

Another important property of a payoff allocation is the \emph{individual rationality}, which requires
\begin{equation}
x_i \geq v(\{i\}) = 0, \ \forall i \in \mathcal{M}^0.
\end{equation}
In order for the Shapley value and the nucleolus to meet this requirement, the value function needs to be monotone \cite{Karaca2019}: $\mathcal{T}_1 \subseteq \mathcal{T}_2 \Rightarrow v(\mathcal{T}_1) \leq v(\mathcal{T}_2)$. However, due to the uncertainty associated with the realization of the customers' schedulable loads, it is possible for some customer's data to drive the retailer's update of wholesale quantity in the suboptimal direction. Therefore, this game is not guaranteed to be monotone. As a result, the Shapley value and the nucleolus can be negative for some players in certain timesteps, as shown in Fig. \ref{fig:sh_t}. However, the case study results presented in Section \ref{sec:case_studies} provide some empirical evidence that these payoff allocations, especially the Shapley value, are likely to be positive in expectation. The proof of this hypothesis is an interesting extension of this work.

\section{Case Studies} \label{sec:case_studies}

Our proposed cooperative game theoretic framework monetizes the customer load data for energy customers and their energy retailer, who is modeled as a cost-based newsvendor. Due to the probabilistic nature of the model, we first apply this framework to an extended time span to evaluate the long-term financial feasibility for the retailer and their customers; then we demonstrate the impact of different methods to determine the wholesale purchase on the retailer's profit distribution.

\subsection{Payoff Allocation for an Extended Time Span} \label{case_largeNo}

\begin{figure} 
\begin{center}
\includegraphics[width=8.4cm]{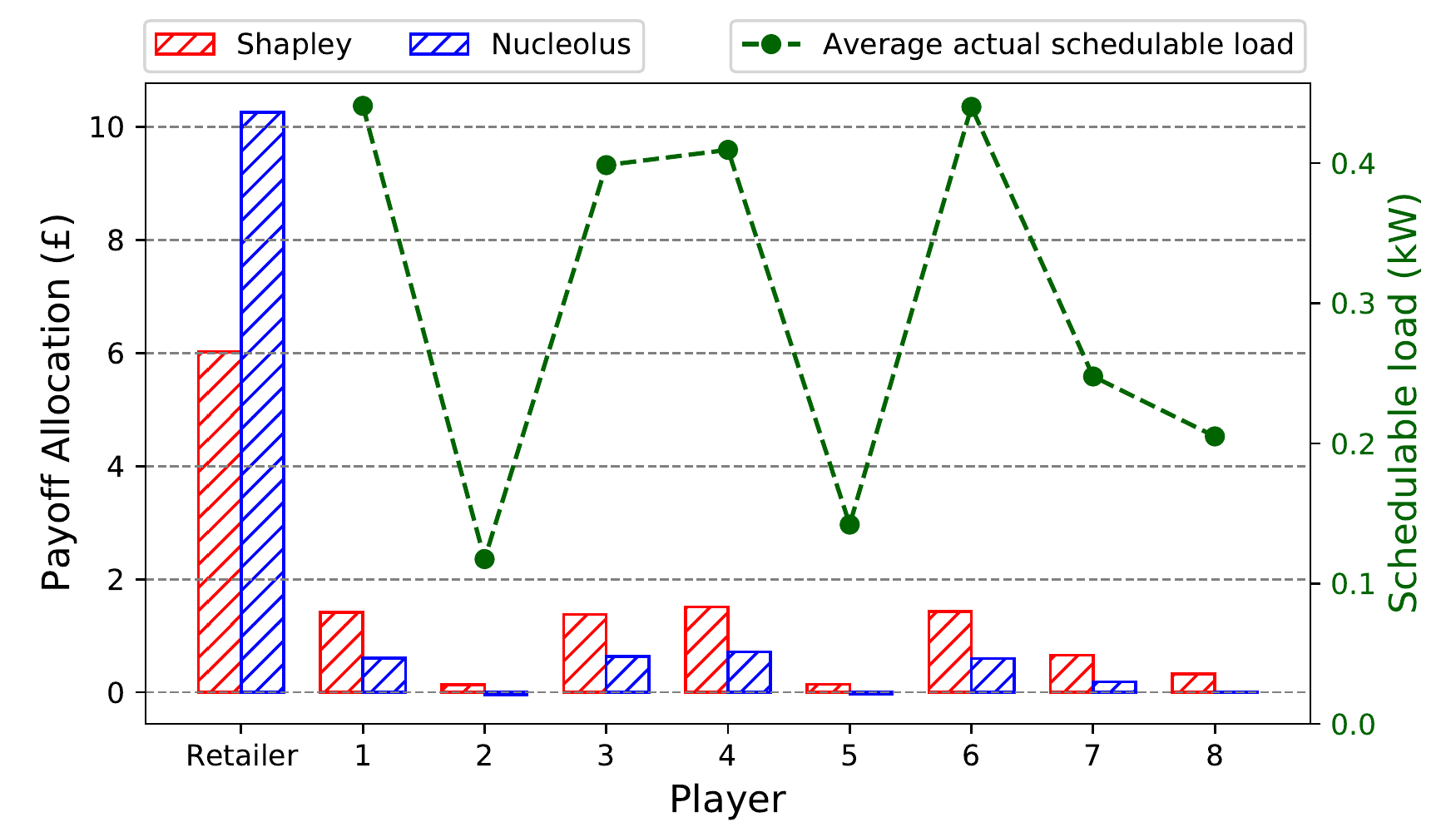}    
\caption{Monthly total profits allocated to the retailer and all the customers (Shapley value vs. nucleolus) and the customers' average schedulable loads.} 
\label{fig:imp_sch_mo}
\end{center}
\end{figure}

Using the same model inputs as previous sections, we extend the time span from 1 day to 1 month. Still using 30 minutes as the time unit for market clearing, we plot the total aggregated payoff allocations by player in Fig. \ref{fig:imp_sch_mo}. Comparing the Shapley value and the nucleolus, it is evident that the Shapley value offers the customers more profit, while the nucleolus favors the retailer. The percentages of the retailer's share of profit are 46\% and 79\%, respectively. This may be explained by the nucleolus' tendency to minimize the dissatisfaction of the retailer because of their essential role in creating value for all the coalitions that have positive values. As a result, some aggregated nucleolus allocations remain negative, even though all the aggregated Shapley allocations are positive. Based on this comparison alone, the Shapley value is the preferred allocation as it maintains a level of incentive for the customers. Another interesting observation is that the magnitude of both payoffs follow similar trends if compared across players, and this trend is consistent with the average schedulable loads. The rationale behind this is that the larger a customer's schedulable load is, the more impactful their schedulable load information is in steering the retailer's wholesale purchase to the optimal amount, and the more additional profit can be created. This chain of influence is captured in the coalition value calculations, which are eventually reflected in the payoff allocations.

\begin{figure} 
\begin{center}
\includegraphics[width=8.4cm]{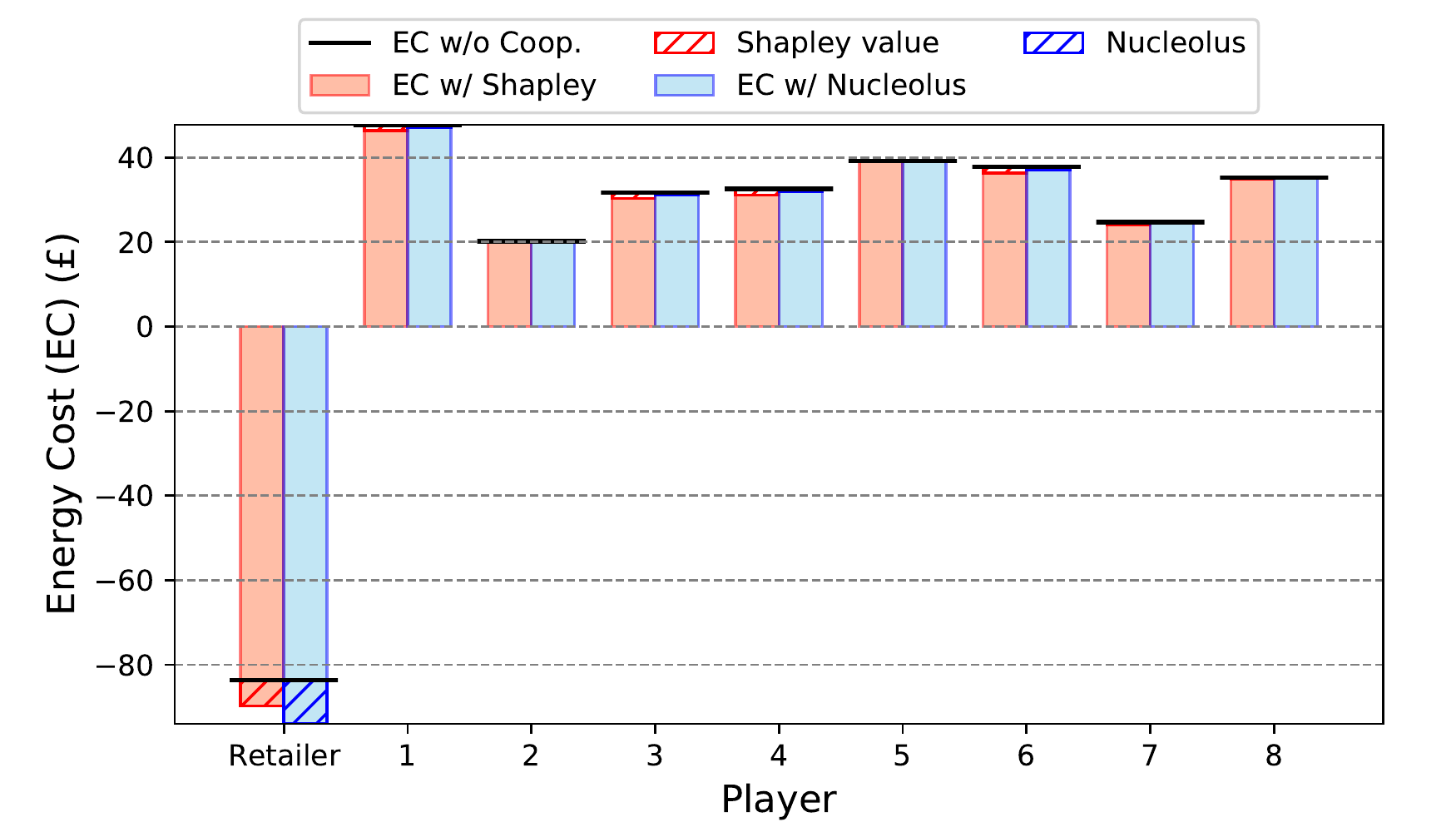}    
\caption{Monthly total energy costs and aggregated profits allocated to the retailer and all the customers.} 
\label{fig:ec_imp_mo}
\end{center}
\end{figure}

Fig. \ref{fig:ec_imp_mo} compares the allocated payoffs with the retailer's profit and the customers' energy costs. The retailer receives an additional 7.2\% and 12.3\% of profit based on the Shapley value and the nucleolus respectively, which is quite significant. The average savings are only 2.6\% and 1.0\% for the customers based on the Shapley value and the nucleolus, respectively. In the case where the payoffs are not enough to offset some customers' privacy concerns, developing a mechanism to readjust the  allocation or to allow customers to reject the offer is an interesting area for future work.

\subsection{Wholesale Quantity and Load Realization}

\begin{figure} 
\begin{center}
\includegraphics[width=8.4cm]{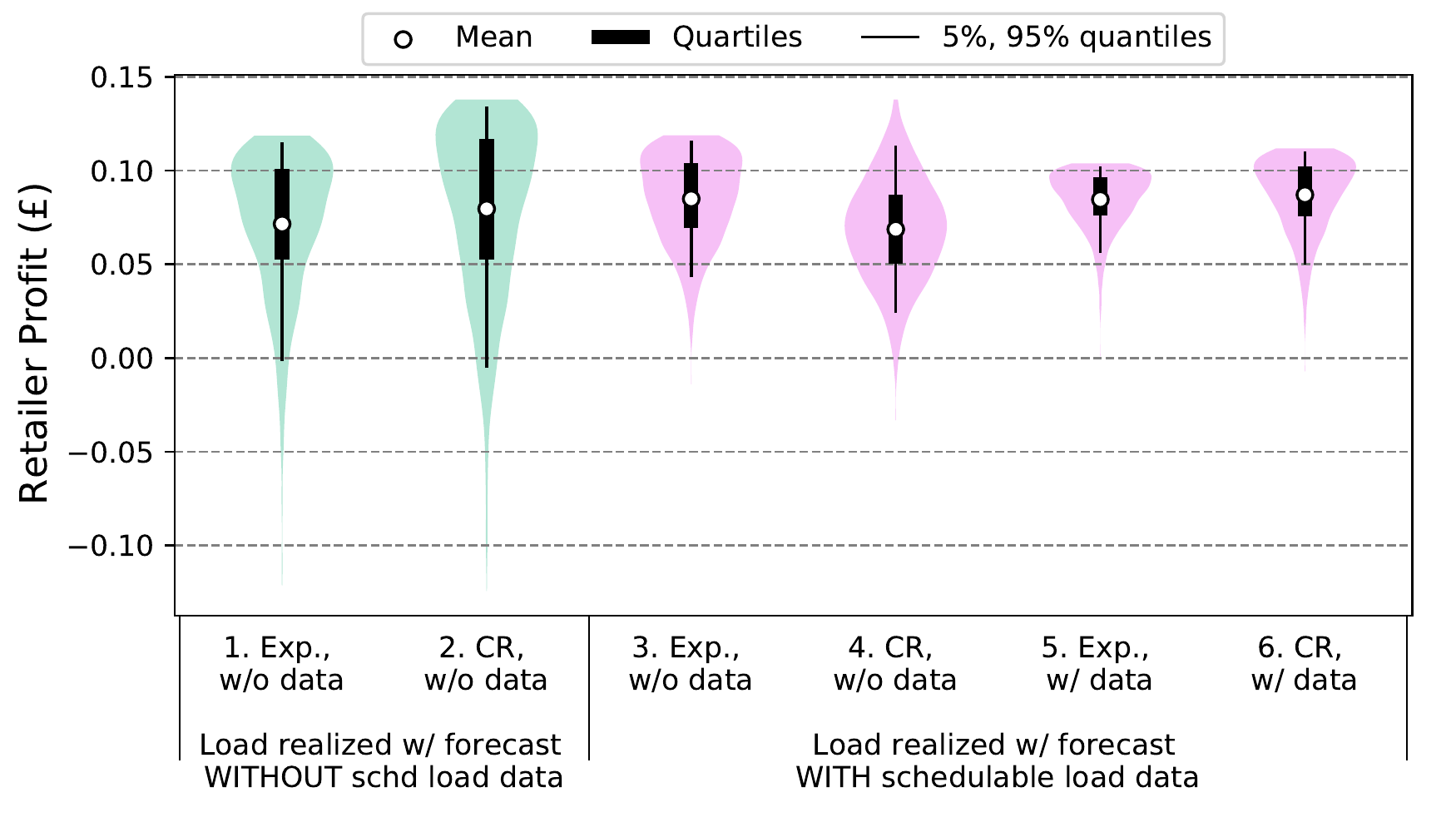}    
\caption{Profit distribution comparison based on forecast and the method to determine the wholesale quantity to purchase.} 
\label{fig:prof_forecast_wholesale}
\end{center}
\end{figure}

The value of a ``narrower'' forecast as shown in Fig. \ref{fig:pdf_cdf_prof} is not only to increase the retailer profit in expectation, but also to reduce the risk of very low profits. The probabilistic distributions of the retailer's profit are plotted in Fig. \ref{fig:prof_forecast_wholesale} differentiated by three factors: 1) whether the customer loads are realized using the probabilistic distribution without schedulable load data (Plots 1, 2) or with schedulable load data (Plots 3, 4, 5, 6); 2) whether the retailer forecasts the load without schedulable load data (Plots 1, 2, 3, 4) or with schedulable load data (Plots 5, 6); 3) whether the method to determine the wholesale purchase is by using the forecast expectation (Plots 1, 3, 5) or the cost ratio (CR) (Plots 2, 4, 6) following \eqref{eq:opt_bid_cl}. 

For the first differentiation factor, comparing Plot 1 with Plot 3, and Plot 2 with Plot 4, we see that by achieving a narrower forecast of the customer loads, the distribution of profit becomes narrower as well. For the second differentiation factor, comparing Plot 4 with Plot 6, we observe additional profits for having a better forecast of the customer loads. For the third differentiation factor, comparing Plot 1 with Plot 2, and Plot 5 with Plot 6, we show an increase in the profit expectation by using the CR as a method to determine the wholesale purchase. However, as the imbalance prices are unknown during the wholesale market stages, computing the CR is difficult in practice. To resolve this issue, a method to forecast the imbalance prices is needed, or the expectation of the forecast can be used despite reduced savings (Plot 5).

\section{Conclusion} \label{sec:conc}

Modeling the energy retailer as a cost-based newsvendor, we demonstrate the effectiveness of determining the wholesale purchase based on probabilistic forecasting of customer loads in maximizing their profit in expectation. Customer schedulable load data can improve load forecasting, which has a positive impact on the retailer's profit. Therefore, to encourage the energy customers to disclose their schedulable load data, we propose a cooperative game theoretic approach for the retailer to allocate a portion of their additional profit to the customers as monetary incentives. Using the Shapley value and the nucleolus as payoff allocation methods in the case studies, we have shown empirically that the retailer is able to retain over 40\% of the profit, and the monthly aggregated Shapley allocations for all customers are positive. Interesting future work includes the investigation of the payoff allocations' computation complexities and their long-term values in expectation, and the incorporation of the imbalance price forecasting in the retailer's cost function.


\bibliographystyle{myIEEEtran}
\bibliography{references.bib} 

\end{document}